\documentclass[11pt]{article}
\usepackage{graphicx}
\usepackage{psfrag}
\usepackage{amsmath}
\usepackage{amssymb}

\newtheorem{theorem}{Theorem}

\newtheorem{lemma}[theorem]{Lemma}
\newtheorem{corollary}[theorem]{Corollary}
\newtheorem{definition}{Definition}
\newtheorem{conjecture}{Conjecture}

\title{Lissajous Knots and Knots with Lissajous Projections}

\author{
Jim Hoste\thanks{Supported by NSF grant 0453284, Claremont Colleges REU Site}\\
Pitzer College
\and
Laura Zirbel\thanks{Supported by NSF grant 0453284, Claremont Colleges REU Site}\\
Loyola Marymount University
}

\begin{document}
\maketitle

\begin{abstract}
Knots in $\mathbb R^3$ which may be parameterized by a single cosine function in each coordinate are called {\it Lissajous} knots. We show that twist knots are Lissajous knots if and only if their Arf invariants are zero. We further prove that all 2-bridge knots and all $(3,q)$-torus knots have Lissajous projections. 
\end{abstract}

\section{Introduction}
\begin{definition}
\label{Lissajous definition} A Lissajous knot $K$in $\mathbb R^3$ is  one that can be represented parametrically by
\begin{eqnarray*}
x(t)&=&\cos(n_x t+\phi_x)\\
y(t)&=&\cos(n_y t+\phi_y)\\
z(t)&=&\cos(n_z t+\phi_z)
\end{eqnarray*}
for integer frequencies $n_x, n_y$, and $n_z$, real phase shifts $\phi_x, \phi_y$, and $\phi_z$, and $0 \le t \le 2 \pi$. 
\end{definition}
In order for the parameterized curve to actually be a knot, that is, to not have any self-intersections, it follows easily that the  three frequencies, $n_x, n_y$ and $n_z$ must be pair-wise relatively prime. Furthermore, there are restrictions on the phase shifts. If we assume that $\phi_z=0$, then  the $x$ and $y$ phase shifts must {\sl not} be of the form
$$\phi_x=k \frac{\pi}{n_z}, \mbox{ or } \phi_y=k \frac{\pi}{n_z}, \mbox{ or } \phi_x=\frac{n_x}{n_y} \phi_y+k\frac{\pi}{n_y}.$$ 
(See \cite{BHJS} or \cite{L}.) Note that by replacing $t$ with $t+c$ we may change the phase shifts, making any one of them zero for example.  

Lissajous knots were first studied in \cite{BHJS} where a few examples were given and a few properties of Lissajous knots were derived. In particular it was shown that a Lissajous knot with all odd frequencies is strongly plus amphicheiral and one with a single even frequency is $2$-periodic.  These symmetries are easy to see. In the first case, consider the orientation reversing homeomorphism $\rho: \mathbb R^3 \to \mathbb R^3$  given by negating each coordinate, or in other words, reflecting through the origin. 
If $K(t)$ is a Lissajous knot with all odd frequencies then
\begin{eqnarray*}
\rho(K(t))&=&(-\cos(n_x t+\phi_x),  -\cos(n_y t+\phi_y), -\cos(n_z t+\phi_z))\\
&=&(\cos(n_x( t+\pi)+\phi_x),  \cos(n_y ( t+\pi)+\phi_y), \cos(n_z( t+\pi)+\phi_z))\\
&=&K(t+\pi).
\end{eqnarray*}
Furthermore, it is not difficult to check that the orientation of $K$ is preserved. Thus $K$ is strongly plus amphicheiral. If one of the frequencies is even, say $n_x$, then a similar argument will show that a 180 degree rotation of $\mathbb R^3$ around the $x$-axis will take $K$ to itself. Thus $K$ is 2-periodic.

Being strongly plus amphicheiral is a fairly restrictive property and among all prime knots to 12 crossings only 3 have this quality~\cite{HTW}. These are 10a103 ($10_{99})$, 10a121 ($10_{123}$), and 12a427.   Knot names are given in both the  Dowker-Thistlethwaite ordering of the Hoste-Thistlethwaite-Weeks table  \cite{HTW} and, in parenthesis,  the Rolfsen~\cite{R} ordering (for knots with 10 or less crossings). Symmetries of the knots in the HTW table were computed using {\it SnapPea} as described in \cite{HTW}. The fact that 10a103 is Lissajous is given in~\cite{L} and to date, it is the only one of these knots that has been shown to be Lissajous.

If $K$ is any knot and $\overline K$ is its mirror image, then $K \sharp \overline K$ is strongly plus amphicheiral. Currently, only three Lissajous knots of this form have been found, namely when $K$ is 3a1 ($3_1$), 5a1 ($5_2$), and 6a3 ($6_1$). A list of Lissajous knots with small crossing number which includes these examples is given in \cite{L}.

In the case where one of the frequencies is even, and the knot is $2$-periodic, actually more is known. In this case, it is shown in both  \cite{L} and \cite{JP} that the knot must have linking number $\pm 1$ with the rotational axis of symmetry. Both the trefoil knot and the figure eight knot for example are $2$-periodic. But neither has such a symmetry where their linking number with the axis of rotation is $\pm 1$. This is because, as we will see shortly,  neither can  be a Lissajous knot! 

The symmetry of a Lissajous knot has a strong effect on its Alexander polynomial. Hartley and Kawauchi~\cite{HK} show that a strongly plus amphicheiral knot must have an Alexander polynomial which is a perfect square. For knots that are $2$-periodic and link their axis of rotation $\pm 1$ times,  Murasugi ~\cite{M} has shown that the Alexander polynomial must be a square modulo 2. From these facts, Lamm shows that many large classes of knots are not Lissajous. (See also \cite{JP}.) For example, if $p$ and $q$ are relatively prime integers larger than 1, then the $(p,q)$-cable of any knot is not Lissajous. In particular, torus knots, iterated torus knots and non-trivial algebraic knots cannot be Lissajous. (Thus the trefoil is not Lissajous.) Lamm also points out that these conditions on the Alexander polynomial also imply that fibered knots with odd genus, and fibered 2-bridge knots are not Lissajous. (Thus the figure eight knot is not Lissajous.)

It is well known that the Arf invariant of a knot may be derived from its Alexander polynomial (see the proof of Lemma~\ref{alexpolyformula}), and in \cite{L} it is shown that if the Alexander polynomial is a square modulo 2 then the Arf invariant is zero. Thus the trefoil and figure eight knots cannot be Lissajous as each has an Arf invariant of 1. The fact that the Arf invariant  of a Lissajous knot is zero was also proven in \cite{BHJS} by a different method.

All of the necessary conditions mentioned so far for a knot to be Lissajous are derived from the symmetry of the knot: either strongly plus amphicheiral, or 2-periodic with linking number $\pm 1$ with the axis of rotation. In general, the implied statements about the Alexander polynomial are strictly weaker and the fact that the Arf invariant must be zero is weaker still.\footnote{ The knot 7a7 ($7_1$ in Rolfsen) has Arf invariant zero but its Alexander polynomial is not a square modulo 2. The knot 8n1 ($8_{20}$ in Rolfsen) has Alexander polynomial which is a square modulo 2 but it is not $p$-periodic  for any $p$. (See Table I in \cite{BZ}). Finally, the knot 8n1 ($8_{20}$ in Rolfsen) has Alexander polynomial a perfect square, but it is not strongly plus amphicheiral.} To date these are the only known necessary conditions for a knot to be Lissajous and it seems unlikely that  these necessary conditions are sufficient.   However, in this paper we show that these conditions are sufficient for the class of twist knots. In fact,  a twist knot is Lissajous if and only if its Arf invariant is zero. For those twist knots that are not Lissajous we know that infinitely many are {\it second order} Lissajous. A second order Lissajous knot is defined as follows.
\begin{definition} A second order Lissajous knot $K$ is  one that can be represented parametrically by
\begin{eqnarray*}
x(t)&=&\cos(n_1 t+\phi_1)\\
y(t)&=&\cos(n_2 t+\phi_2)\\
z(t)&=&\cos(n_3 t+\phi_3)+\cos(n_4 t+\phi_4)
\end{eqnarray*}
for integer frequencies $n_1, n_2, n_3$, and $n_4$, real phase shifts $\phi_1, \phi_2, \phi_3$, and $\phi_4$, and $0 \le t \le 2 \pi$.
\end{definition}
We believe the following conjecture to be true, and hope to include this result in a forthcoming paper.
\begin{conjecture} Every twist knot with Arf invariant 1 is second order Lissajous.
\label{Arf1conjecture}
\end{conjecture}
Note that a second order Lissajous knot has a {\it Lissajous projection} (in the $x$-$y$ plane). By employing large enough  frequencies it seems reasonable to make the following conjecture.
\begin{conjecture} Every knot has a Lissajous projection.
\label{LPconjecture}
\end{conjecture}
While not explicitly stated as a conjecture in \cite{BHJS}, this belief is evident in their paper and we echo their sentiment here.

If Conjecture~\ref{LPconjecture} is true, then the height function $z(t)$ can be estimated by a finite sum of cosine functions. Whether just two (or one) cosine functions would then suffice in all cases, as we believe it does for all twist knots, is an interesting question.

In Section~\ref{twist knots}  we show that a twist knot is Lissajous if and only if  its Arf invariant is zero. Next, in Section~\ref{LP} we give a general result that can be used to show that a knot has a Lissajous projection. We use this to show that all 2-bridge knots have Lissajous projections.  Finally, in the last section we continue with the techniques in Section~\ref{LP} to show that all $(3,q)$-torus knots have Lissajous projections.

{\bf Acknowledgement} This research was carried out during the summer of 2005 as part of a {\sl Research Experiences for Undergraduates} program administered by the Claremont Colleges with funding from the National Science Foundation. The authors wish to thank the Claremont Colleges for their support and hospitality as well as both the Claremont Colleges and the National Science Foundation for funding the REU site.

{\bf Postscript} After completing this paper we learned that Lamm has proven Conjecture~\ref{LPconjecture}.

\section{Twist Knots}
\label{twist knots}
\begin{figure}[b]
    \begin{center}
    \leavevmode
    \scalebox{.70}{\includegraphics{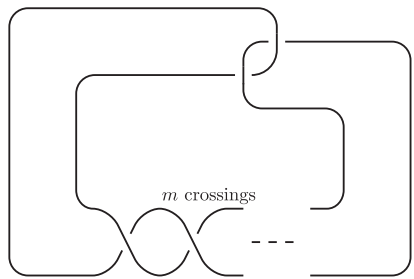}}
    \end{center}
\caption{The twist knot $K_m$.}
\label{twistKnot}
\end{figure}
The twist knot $K_m$ is defined as in Figure~\ref{twistKnot}. There are $m$ crossings in the ``twisty''  part and then two more crossings in the left-handed ``clasp.''  If $m>0$ then the given diagram is alternating (and reduced) and thus depicts a knot with crossing number equal to $m+2$. If $m<0$, then the crossings in the ``twisty'' part have the opposite handedness. For example, $K_2$ is the figure eight knot, $K_1$ is the left-handed trefoil, and $K_{-1}$ is the unknot. It is not difficult to show that $K_m$ and $K_{-m-1}$ are mirror images of each other. We will consider two knots the same if there is a homeomorphism of $\mathbb R^3$ to itself taking one knot to the other and therefore will not distinguish between a knot and its mirror image. Thus, it suffices to consider only the case where $m=2n$ is even.

 We begin with the following Lemma. 
\begin{lemma} The Alexander polynomial  and the Arf invariant of the twist knot $K_{2n}$ are
\label{alexpolyformula} 
\begin{eqnarray*}
\Delta_{K_{2n}}(t)&=& n-(2n+1)t+n t^2\\
\mbox{Arf}\ (K_{2n})&\equiv&n \mbox{ mod}\ 2.
\end{eqnarray*}
\end{lemma}
\noindent{\bf Proof:} The first part is well known. See \cite{R} for example. 

For any knot $J$, the Arf invariant may be derived from the Alexander polynomial as follows.
$$\mbox{Arf}(J)=\begin{cases}
0&\text{if $ \Delta_J(-1)\equiv \pm 1$  mod 8,}\\
1&\text{if $ \Delta_J(-1)\equiv \pm 5$ mod 8.}
\end{cases}
$$
In our case, $\Delta_{K_{2n}}(-1)=4 n+1$ which is 1 mod 8 if $n$ is even and 5 mod 8 if $n$ is odd.
\hfill $\square$

\begin{figure}[h]
\psfrag{n}{$n$}
    \begin{center}
    \leavevmode
    \scalebox{.75}{\includegraphics{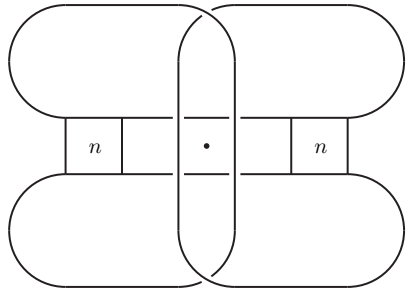}}
    \end{center}
\caption{$K_m$ is 2-periodic and links its axis of rotation once iff $m\equiv 0$ mod 4.}
\label{symmetry}
\end{figure}

From Lemma~\ref{alexpolyformula} we see that the Alexander polynomial of a twist knot is never a perfect square and hence twist knots are never strongly plus amphicheiral. Thus none can be   Lissajous knots with all odd frequencies. It remains to consider the possibility that one frequency is even. The following lemma narrows the scope of our investigation.
\begin{lemma} For twist knots $K_m$, the following are equivalent:
\label{equivalent conditions for twist knots}
\begin{enumerate}
\item  $K_m$ is 2-periodic and has linking number $\pm 1 $ with its rotational axis of symmetry.
\item The Alexander polynomial $\Delta_{K_m}(t)$ is a square modulo 2.
\item The Arf invariant of $K_m$ is zero.
\end{enumerate}
\end{lemma}
\noindent {\bf Proof:} For any knot, $1 \Rightarrow 2  \Rightarrow 3$. Thus it remains to prove $3  \Rightarrow 1$. As mentioned already, we may assume that $m=2n$ is an even integer.
Suppose the Arf invariant of $K_{2n}$ is zero.   Thus Lemma~\ref{alexpolyformula} implies that $n$ must be even. Figure~\ref{symmetry} now shows that $K_m$ is 2-periodic. In the figure, the box labeled $n$ indicates a total of $n$ crossings. The axis of rotation is perpendicular to the plane of the diagram and appears as a heavy dot in the center of the figure. The reader should check that the linking number with the axis is $\pm 1$ if $n$ is even and $\pm 3$ if $n$ is odd.
\hfill $\square$

\begin{figure}
    \begin{center}
    \leavevmode
    \scalebox{.75}{\includegraphics{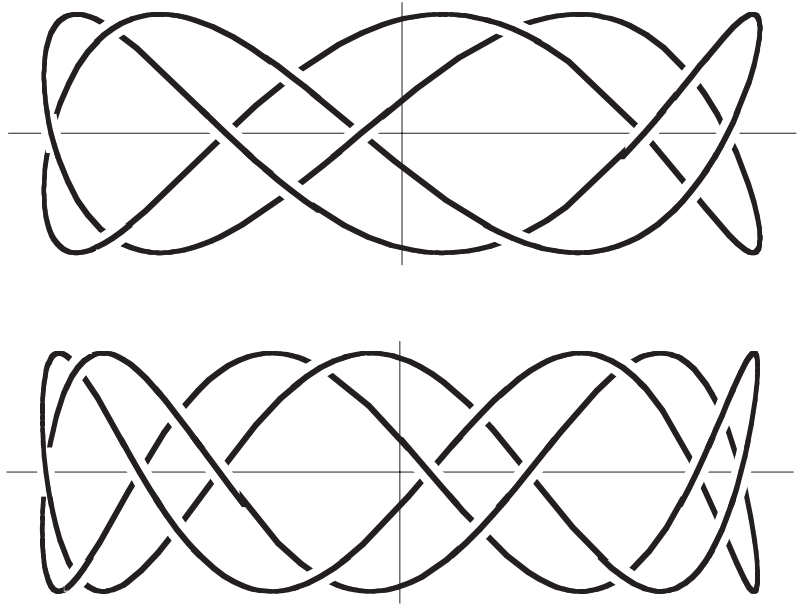}}
    \end{center}
\caption{Lissajous knots with $n_x=2, n_y=2m+1, n_z=6m+7, \phi_x=0, \phi_y=0.5$ and $\phi_z=(6m+7-3 \pi)/(2(2m+1))$ for $m=2$ and 3.}
\label{m=2 or 3}
\end{figure}

Because of Lemma~\ref{equivalent conditions for twist knots} it is quite natural to ask whether or not twist knots with Arf invariant equal to zero are Lissajous. By exploring Lissajous knots with small frequencies one quickly finds that if $n_x=2$ and $n_y$ is odd, then the Lissajous knot has a projection which {\sl could} be the projection of a twist knot if all of the crossings are arranged in the appropriate way. For example, consider the two Lissajous knots shown in Figure~\ref{m=2 or 3} where $n_x=2, \phi_x=0, n_y=5$ or 7 and $\phi_y=0.5$. Arbitrarily choosing $\phi_y=0.5$, we are led to the following question: If $n_x=2, \phi_x=0, n_y=2m+1$ and $\phi_y=0.5$, can we choose $n_z$ and $\phi_z$ so as to produce a twist knot? In order to study this question we begin by investigating the double points of the $xy$-projection of a general Lissajous knot.
\begin{lemma} Let $K(t)$ be the Lissajous knot parameterized  as in Definition~\ref{Lissajous definition}. Suppose $t=t_1$ and $t=t_2$ give  a double point in the projection to the  $xy$-plane. Then the set of such pairs $(t_1, t_2)$ fall into two types.
{\renewcommand{\arraystretch} {1.8} 
\begin{tabular}{lrl}
Type I:& $(t_1, t_2)=$&$((-\frac{k}{n_x}+\frac{j}{n_y})\pi-\frac{\phi_y}{n_y},(\frac{k}{n_x}+\frac{j}{n_y})\pi-\frac{\phi_y}{n_y})$\\
&& for $1\le k\le n_x-1$ and $\lceil \frac{n_y}{n_x}k+\frac{\phi_y}{\pi}\rceil \le j \le \lfloor 2 n_y-\frac{n_y}{n_x}k+\frac{\phi_y}{\pi}\rfloor$.\\

Type II:& $(t_1, t_2)=$&$((-\frac{k}{n_y}+\frac{j}{n_x})\pi-\frac{\phi_x}{n_x},(\frac{k}{n_y}+\frac{j}{n_x})\pi-\frac{\phi_x}{n_x})$\\
&& for $1\le k\le n_y-1$ and $\lceil \frac{n_x}{n_y}k+\frac{\phi_x}{\pi}\rceil \le j \le \lfloor 2 n_x-\frac{n_x}{n_y}k+\frac{\phi_x}{\pi}\rfloor$.\\
\end{tabular}
}
Moreover, there are $n_xn_y-n_y$ doubles points of Type I, and $n_xn_y-n_x$ double points of Type II. 
\end{lemma}
\noindent{\bf Proof:} We must find all pairs  of times $(t_1, t_2)$ such that $x(t_1)=x(t_2)$ and $y(t_1)=y(t_2)$. Using the fact that $\cos u=\cos v$ if and only if $u\pm v$ is an integer multiple of $2 \pi$, one is led to a number of cases that must be analyzed. See \cite{BHJS} for further details.
\hfill $\square$

Let $C^I_{k,j}$ and $C^{II}_{k,j}$ denote the crossings that occur at the double points of Type I and II respectively. By $x(C^I_{k,j})$ or $y(C^I_{k,j})$ we will mean the $x$ and $y$ values at that crossing.

Suppose now that $n_x=2, \phi_x=0, n_y=2m+1$ and $\phi_y=0.5$. Then there are $2m+1$ double points of Type I. For each of these, $k=1$ and $m+1\le j\le 3m+1$. The time pairs are
$$(t_1, t_2)_I=\left ( 
\frac{(2j-2m-1)\pi-1}{2(2 m+1)}, 
\frac{(2j+2m+1)\pi-1}{2(2 m+1)} 
\right ).$$
It is easy to show that $y(C^I_{1,j})=0$ for all $j$ and hence each Type~I crossing lies on the $x$-axis.
Furthermore, it is not difficult to show that no Type~I crossing can occur at the origin and the number of Type~I double points that fall on the positive $x$-axis is 
$2\left \lceil  \frac{m}{2} \right \rceil$.

If we now look at the Type~II double points we see that there are $4m$ of these given by the time pairs
$$(t_1, t_2)_{II}=\left (
\frac{(2m+1)j-2k}{2(2m+1)}\pi,
\frac{(2m+1)j+2k}{2(2m+1)}\pi
\right)$$
where $1 \le k \le m$ if $j=1$ or $3$ and $1 \le k \le 2 m$ if $j=2$. If we compute the value of $y$ at these crossings we obtain
\begin{eqnarray*}
y(C^{II}_{k,1})&=&(-1)^{k+m+1} \sin\frac{1}{2}\\
y(C^{II}_{k,2})&=&(-1)^{k+1} \cos\frac{1}{2}\\
y(C^{II}_{k,3})&=&-y(C^{II}_{k,1}).
\end{eqnarray*}
By computing the $x$ value of the Type~II crossings we obtain
\begin{eqnarray*}
x(C^{II}_{k,1})&=&x(C^{II}_{k,3})\\
x(C^{II}_{k,2})&=&x(C^{II}_{2 m+1-k, 2})\\
\end{eqnarray*}
for all $k$. Furthermore, the number of Type~II crossings $C^{II}_{k,j}$ with $x>0$ is $2 m$. Of these, $2\left \lfloor  \frac{m}{2} \right \rfloor$ have $j=2$ while $2\left \lceil  \frac{m}{2} \right \rceil$ have $j=1$ or $3$.

Finally, since $n_x=2$, the knot is rotationally symmetric around the $x$-axis. Thus the signs of $C^{II}_{k,1}$ and $C^{II}_{k,3}$ are the same as are the signs of $C^{II}_{k,2}$ and $C^{II}_{2m+1-k, 2}$.

We now want to understand how the choice of $n_z$ and $\phi_z$ determines the crossing at each of the double points in the $xy$-projection. Note first that  increasing the phase shift by $\pi$ will reflect the knot though the $xy$-plane. Thus we may assume that $0\le \phi_z \le \pi$.

Suppose that at a double point of either type, $z(t_1)=z(t_2)$. This occurs if and only if the sum or difference of $n_z t_1+\phi_z$ and $n_z t_2+\phi_z$ is an integer multiple of $2 \pi$. If we examine the Type~II double points first, we find that all of the crossings $C^{II}_{k,2}$ become singular precisely when $\phi_z=0$ or $\pi$. All the other Type~II crossings become singular precisely when $\phi_z=\frac{\pi}{2}$. Furthermore, both of these statements are independent of $n_z$. In particular, once $n_z$ is chosen, we are free to vary $\phi_z$ between $0$ and $\frac{\pi}{2}$ without changing any of the Type~II crossings! We now turn our attention to the Type~I crossings. Since $n_z$ is odd, it is not hard to show that the difference of $n_z t_1+\phi_z$ and $n_z t_2+\phi_z$ cannot be an even multiple of $\pi$. However, their sum  is an even multiple of $\pi$ if and only if
$$\phi_z=\frac{n_z}{2(2m+1)}-\frac{jn_z \pi}{2m+1}+i \pi$$
for some integer $i$. Therefore, there are $2 m+1$ values of $\phi_z$ between $0$ and $\pi$ where all the Type~I crossings become singular. These values, together with $\frac{\pi}{2}$, give $2m+2$ phases that must be avoided. The locations of these singular phases depend on $n_z$ but the number of them depend only on $m$. Thus, for a given choice of $m$ there are $2 m+3$ intervals where we may safely choose $\phi_z$ and varying $\phi_z$ within each interval will not cause any of the crossings to become singular. Hence there are at most $2 m+3$ different Lissajous knot types with $n_x=2, n_y=2 m+1, \phi_x=0$ and $\phi_y=\frac{1}{2}$.

\begin{figure}[h]
\psfrag{a}{a}
\psfrag{b}{b}
\psfrag{c}{c}
\psfrag{d}{d}
    \begin{center}
    \leavevmode
    \scalebox{1.0}{\includegraphics{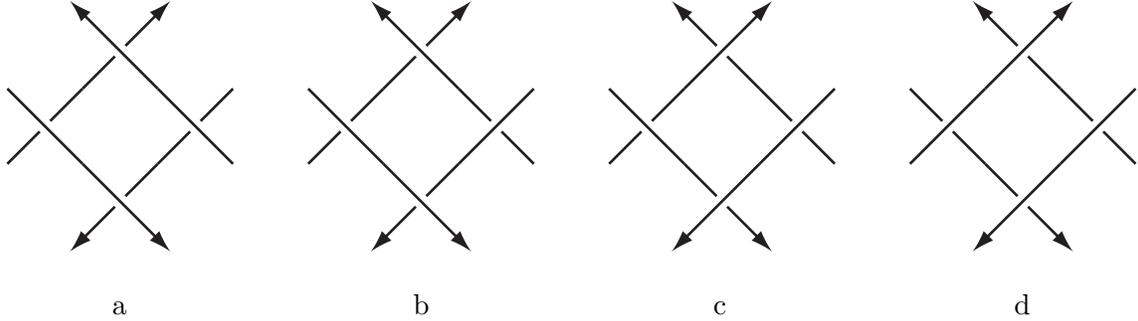}}
    \end{center}
\caption{a: left band crossing; b: clasp, $m$ even; c: clasp, $m$ odd; d: right band crossing.}
\label{band crossing}
\end{figure}

At this point it is a simple matter to program a computer to search through a large number of pairwise relatively prime frequencies $2, 2m+1, n_z$ and for each of these construct the $2m+3$ possible Lissajous knots. Looking for twist knots in this way led to the following theorem.
\begin{theorem} 
\label{arf0iffLissajous} A twist knot is Lissajous if and only if its Arf invariant is zero. In particular, the Lissajous knot
\begin{eqnarray*}
x(t)&=&\cos(2 t) \\
y(t)&=&\cos((2m+1) t+0.5)\\
z(t)&=&\cos((6m+7) t+(6m+7-3 \pi)/(4 m+2))
\end{eqnarray*}
for $m\ge 0$ is equivalent to $K_{2 m}$ if $m$ is even and $K_{-2(m+1)}$ if $m$ is odd.
\end{theorem}
\noindent{\bf Proof:} If $m=0$ or $m=1$ it is easy to verify the result. Thus, in everything that follows  we assume that $m>1$.  Diagrams for $m=2$ and 3 are shown in Figure~\ref{m=2 or 3}. A nice pattern exists that we will show persists for all $m$. 
Consider the checkerboard coloring of the diagram and assume the unbounded region has been colored white. There are exactly $m$ white regions each  having 4 edges and 4 vertices. We call these {\it diamonds}.  Each diamond lies on the $x$-axis and contains  two Type~I crossings and two Type~II crossings.  There are a total of $m$ diamonds and when $m$ is even, half lie on each side of the $y$-axis. When $m$ is odd, $(m+1)/2$ lie in the right half-plane and $(m-1)/2$ lie in the left half-plane. We claim that each diamond must appear as one of the four possibilities shown  in Figure~\ref{band crossing}, which we call {\it left band crossings}, {\it clasps}, or {\it right band crossings}. 

In Figure~\ref{m=2 or 3} for example, notice that as we move from left to right the diamonds form a (possibly empty) string of left band crossings, then a clasp, and then a string of right band crossings. Thus all of the right band crossings ``unravel'' and the clasp of the twist knot is formed at the diamond that is of type b) or c) in Figure~\ref{band crossing}. From this it clearly follows that the knot is indeed a twist knot. 

It suffices to prove the following four claims:
\begin{enumerate}
\item If $m$ is even, then all Type~I crossings in the right half-plane are left-handed and all Type~I crossings in the left half-plane, {\sl except} the one closest to the $y$-axis, are right-handed.
\item If $m$ is odd, then all Type~I crossings in the left half-plane are right-handed and all Type~I crossings in the right half-plane, {\sl except} the one closest to the $y$-axis, are left-handed.
\item Type II crossings with $k=2$ are left-handed if and only if they lie in the right half-plane.
\item Type II crossings with $k=1$ or $3$ are right-handed if and only if they lie in the right half-plane.
\end{enumerate}
From these claims it is easy to establish
\begin{enumerate}
\setcounter{enumi}{4}
\item If $m$ is even, the first diamond to the left of the $y$-axis is a clasp as shown in Figure~\ref{band crossing}b. If $m$ is odd, the first diamond to the right of the $y$-axis is a clasp as shown in Figure~\ref{band crossing}c. In either case, all diamonds to the right of this clasp are right band crossings and all diamonds to its left are left band crossings.
\end{enumerate}
Once statements 1--5 are proven we obtain the desired result as follows. Each of the left band crossings contributes one full left-handed twist to the {\it band} of the twist knot and each of the Type~II crossings with $k=2$ in the left half-plane contributes  half a left-handed twist. Notice that since the edges of the band are oriented oppositely, these are actually right-handed crossings. Finally, the left-most Type~I crossing is right-handed and contributes half a left-handed twist to the band. Thus if $m$ is odd for example, all of these crossings add up to $-2m-1$ right half-twists in the band. The clasp of the twist knot is right-handed. Thus this knot is the mirror image of $K_{2m+1}$ which is equivalent to $K_{-2(m+1)}$.

In order to establish statements 1--4 we study the signs of the crossings. Suppose $(t_1, t_2)$ is a pair of times that gives a double point, and let $$v_i=(x'(t_i), y'(t_i), 0)$$ be the tangent vector to the $xy$-projection at $t=t_i$. The cross product
$$v_1 \times v_2=\left (0,0, x'(t_1)y'(t_2)-x'(t_2)y'(t_1)\right )$$
 is now a vector perpendicular to the $xy$-plane. The following lemma is proven in \cite{L}.
 \begin{lemma} \label{sign of crossing}The sign of the crossing corresponding to the times $(t_1,t_2)$ is
$$\text{sign}\left [(x'(t_1)y'(t_2)-x'(t_2)y'(t_1))(z(t_1)-z(t_2))\right].$$
\end{lemma}
Using Lemma~\ref{sign of crossing},  somewhat lengthy, but nevertheless straightforward, calculations show that statements 1--4 are correct. We illustrate the proof of statement 3 and leave the others to the reader.

Considering a Type~II crossing with $j=2$, we obtain
$$x'(t_2)=-x'(t_1)=2 \sin\frac{2 k \pi}{2 m+1} \mbox{ and}$$
$$y'(t_2)=y'(t_1)=(-1)^k(2 m+1)\sin\frac{1}{2}.$$
Hence the sign of  $x'(t_1)y'(t_2)-x'(t_2)y'(t_1)$ is the same as the sign of
$$(-1)^k\sin\frac{2 k \pi}{2 m+1}.$$

Next, we consider $z(t_1)-z(t_2)$. As mentioned earlier, varying $\phi_z$ between $0$ and $\pi/2$ will not change the Type~II crossings with $j=2$. Thus we may replace $z(t)$ with $Z(t)=\cos((6 m+7)t+\pi/4)$. We obtain,
\begin{eqnarray*}
Z(t_1)&=&(-1)^{k+1}\frac{1}{\sqrt{2}}(\cos\frac{4 k \pi}{2 m+1}+\sin\frac{4 k \pi}{2 m+1})\\
Z(t_2)&=&(-1)^{k+1}\frac{1}{\sqrt{2}}(\cos\frac{4 k \pi}{2 m+1}-\sin\frac{4 k \pi}{2 m+1})
\end{eqnarray*}
with the sign of $Z(t_1)-Z(t_2)$ equal to $(-1)^{k+1}\mbox{sign}\left [\sin\frac{4 k \pi}{2 m+1}\right ]$.

Finally, using Lemma~\ref{sign of crossing}, the sign of $C^{II}_{k,2}$ is
\begin{eqnarray*}
\mbox{sign}\left [C^{II}_{k,2}\right]&=&\mbox{sign}\left[(-1)^k\sin\frac{2 k \pi}{2 m+1} (-1)^{k+1}\sin\frac{4 k \pi}{2 m+1}\right]\\
&=&-\mbox{sign}\left[\sin\frac{2 k \pi}{2 m+1} \sin\frac{4 k \pi}{2 m+1}\right]\\
&=&-\mbox{sign}\left[\cos\frac{2 k \pi}{2 m+1}\right]\\
&=&-\mbox{sign}\left[ x(C^{II}_{k,2})\right].
\end{eqnarray*}
\hfill $\square$
  
\section{Knots With Lissajous Projections}
\label{LP}
Theorem~\ref{arf0iffLissajous} shows that every twist knot with Arf invariant equal to zero has a Lissajous projection with $n_x=2$ and $n_y$ some odd integer. Our first goal in this section is to show that all 2-bridge knots have Lissajous projections.  To do this we first introduce the concepts of a {\it Lissajous arc} and a {\it Lissajous braid}.

\begin{figure}
    \begin{center}
    \leavevmode
    \scalebox{.50}{\includegraphics{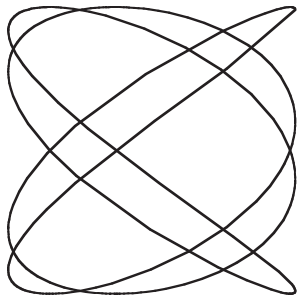}}
    \end{center}
\caption{A Lissajous projection with $n_x=4, n_y=3, \phi_x=0$, and $\phi_y=0.2$.}
\label{Lissajous Projection}
\end{figure}

\begin{figure}
    \begin{center}
    \leavevmode
    \scalebox{.50}{\includegraphics{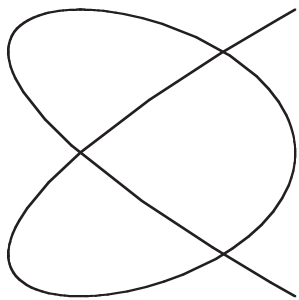}}
    \end{center}
\caption{A Lissajous arc with $n_x=4$ and $ n_y=3$. }
\label{Lissajous Arc}
\end{figure}

Figure~\ref{Lissajous Projection} illustrates the Lissajous projection with frequencies $n_x=4$ and $n_y=3$ and phase shifts $\phi_x=0$ and $\phi_y=0.2$. In Figure~\ref{Lissajous Arc} we have set both phase shifts to zero. Notice that the first figure appears like the projection of a long thin twisted disk or band whose ``core'' is the arc of the second figure. As the phase shift $\phi_y$ is increased from  zero, the arc in the second figure expands to the band of the first figure. A double point in Figure~\ref{Lissajous Arc} becomes four crossings in Figure~\ref{Lissajous Projection}. But additional crossings are also introduced because the band does not lay flat in the plane of the projection. Instead, think of the band as twisted, with these additional crossings located at the extreme points of the arc in both the $x$ and $y$ directions. These observations motivate the following two definitions.
\begin{definition}
A Lissajous arc with frequencies $n_x$ and $n_y$ is the  plane figure parameterized by
$$x(t)=\cos (n_x t) \text{\quad and \quad} y(t)=\cos (n_y t) $$
for $0 \le t \le \pi$, where $n_x$ and $n_y$ are  relatively prime positive integers.
\end{definition}
Notice that if small neighborhoods of all the maxima and minima in  the  $y$ direction are removed from a Lissajous arc, we are left with the projection of a braid on $n_y$ strings. Suppose $\alpha$ is an $n$-string braid expressed as a product of the standard braid generators $\sigma_1, \sigma_2, \dots, \sigma_{n-1}$ and their inverses. Here $\sigma_i$ denotes the braid where string $i+1$ passes over string $i$ as one travels down the braid. Let $\rho(\alpha)$ denote the obvious projection of the braid obtained by replacing $\sigma_i^{\pm 1}$ with a double point which we will denote as $s_i$.
\begin{definition} A Lissajous braid with frequencies $n_x$ and $n_y$ is any braid $\alpha$ on $n_y$ strings whose projection $\rho(\alpha)$ is obtained from a Lissajous arc of the same frequencies by removing small neighborhoods of all the minima and maxima in the $y$ direction.
\end{definition}
The following characterization of Lissajous braids is straightforward and is left to the reader.
\pagebreak
\begin{lemma}
Let $\alpha$ be a Lissajous braid with frequencies $n_x$ and $n_y$. Define $s_\text{even}$ and $s_\text{odd}$ as 
\begin{eqnarray*}
s_\text{even}&=&s_2s_4\dots s_E\\
s_\text{odd}&=&s_1s_3\dots s_O
\end{eqnarray*}
where $E$ and $O$ are the largest even and odd integers less than $n_y$. Then
$$\rho(\alpha)=\left \{ \begin{array}{ll} 
s_\text{even}(s_\text{odd} s_\text{even})^{\frac{n_x}{2}-1} & \mbox{if  $n_x$ is even and $n_y$ is odd} ;\\
(s_\text{odd} s_\text{even})^\frac{n_x-1}{2} & \mbox{if  $n_x$ is odd and $n_y$ is even} ;\\
(s_\text{even} s_\text{odd})^\frac{n_x-1}{2} & \mbox{if  $n_x$ is odd and $n_y$ is odd} .\\
\end{array}\right.$$
The number of crossings in the braid (or arc of the same frequencies) is $(n_x-1)(n_y-1)/2$. The number of crossings in a Lissajous projection with these frequencies is $2n_xn_y-n_x-n_y$.
\end{lemma}
Thus a Lissajous braid has a projection which is an alternating product of $s_\text{even}$ and $s_\text{odd}$. However not all such products are possible. For example, the projection can never begin and end with  $s_\text{odd}$ and the fact that $n_x$ and $n_y$ are relatively prime impose further restrictions as well.

The following two lemmas will play a key part in showing that every 2-bridge knot has a Lissajous projection.
\begin{lemma}
Suppose $\alpha$ is any $n$-string braid with $n \ge 3$ of the form $\alpha=a \sigma_i^{\epsilon_1}\sigma_i^{\epsilon_2}b$ where $\epsilon_i=\pm 1$. Thus
$\rho(\alpha)=\rho(a) s_i^2 \rho(b)$. Then there exist equivalent braids $\alpha^\prime$ and $\alpha^{\prime \prime}$ such that
$$\rho(\alpha^\prime)=\rho(a)(s_{i-1}s_i)^3 \rho(b) \text{ \quad and \quad } \rho(\alpha^{\prime \prime})=\rho(a)(s_is_{i+1})^3 \rho(b).$$
\label{first replacement lemma}
\end{lemma}
{\bf Proof:}
By using  two Type I and one Type III Reidemeister moves we may change $\alpha$ to
$$\alpha^\prime=a \sigma_{i-1}^{-1}\sigma_i^{-1} \sigma_{i-1}^{\epsilon_1}\sigma_i \sigma_{i-1}\sigma_i^{\epsilon_2} b.$$
We now have    $ \rho(\alpha^\prime)=\rho(a)(s_{i-1}s_i)^3 \rho(b)$.
A similar argument gives $\alpha^{\prime \prime}$.
\hfill $\square$

\begin{lemma}
Any 3-string braid $\alpha$ is equivalent to a braid whose projection is an alternating product of $s_1$ and $s_2$. 
\label{3 string braids are nearly Lissajous}
\end{lemma}
{\bf Proof:} We proceed by induction on the length (that is, the number of crossings) of $\alpha$. If $\alpha$ has length one, then $\rho(\alpha)$ equals either $s_1$ or $s_2$ and we are done.

There are now two cases to consider: either $\alpha=\sigma_1^{\pm 1}\beta$ or $\alpha=\sigma_2^{\pm 1}\beta$, where $\beta$ is a shorter braid.   
 
Suppose $\alpha=\sigma_2^{\pm 1}\beta$. By induction, we may assume that $\beta$ is equivalent to a braid, which we continue to call $\beta$, whose projection alternates between $s_1$ and $s_2$.  If $\rho(\beta)$ begins with $s_1$, then we are done, because $\rho(\alpha)$ now alternates between $s_2$ and $s_1$. If $\rho(\beta)$ begins with $s_2$, then $\rho(\alpha)$ begins with  $s^2_2$, which we can replace with  $(s_1s_2)^3$ by Lemma~\ref{first replacement lemma}.   

A similar argument works if $\alpha=\sigma_1^{\pm 1}\beta$. 
\hfill $\square$

\begin{definition}
The closure of a Lissajous arc is the knot projection obtained by connecting the ends of a Lissajous arc with an arc whose interior is disjoint from the Lissajous arc and which introduces no additional crossings. 
\end{definition}
The significance of this definition will be borne out in the following theorem.
\begin{theorem}
\label{closure of LA implies LP}
Suppose the knot $K$ has a projection which is the  closure of a Lissajous arc with frequencies $n_x$ and $n_y$. Then $K$  has a Lissajous projection with the same frequencies. 
\end{theorem}
{\bf Proof:} We demonstrate the idea of the proof with the following example.

\begin{figure}
    \begin{center}
    \leavevmode
    \scalebox{.90}{\includegraphics{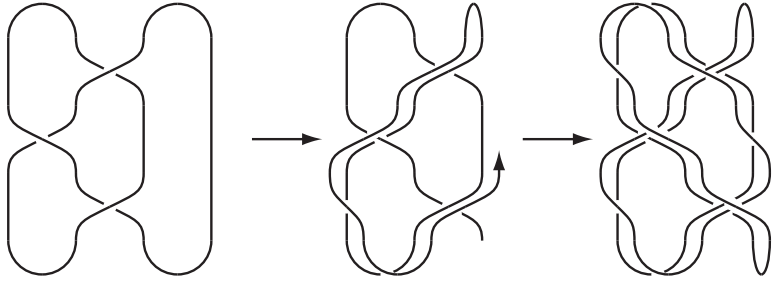}}
    \end{center}
\caption{Changing the closure of a Lissajous arc to a Lissajous projection.}
\label{trefoil2Lissajous}
\end{figure}

Consider the trefoil knot as shown on the left in Figure~\ref{trefoil2Lissajous}. This projection is  the closure of a Lissajous arc. 
The next step is illustrated in the middle and right-hand diagrams in Figure~\ref{trefoil2Lissajous}. Imagine picking up the right-hand strand of the diagram (this is the arc that was used to close the Lissajous arc) and 
then laying it back down atop the figure, retracing the Lissajous arc. As the Lissajous arc is retraced we also introduce crossings as needed to achieve a Lissajous projection. The additional crossings are located at each of the extreme points in both the $x$ and $y$ directions. Note that the final Lissajous projection is equivalent to the one shown in Figure~\ref{Lissajous Projection}.
\hfill $\square$

We now have everything in place to prove the following theorem.
\begin{theorem}
\label{2 bridge implies LP}
Every 2-bridge knot has a Lissajous projection with $n_x=3$.
\end{theorem}
{\bf Proof:} Every 2-bridge knot has a diagram which is the plat closure of a 4-string braid that uses only $\sigma_1, \sigma_2$ and their inverses. We may think of this braid as a 3-string braid $\alpha$ that has had a trivial fourth string added along the right side. By Lemma~\ref{3 string braids are nearly Lissajous} we may alter $\alpha$ so that its projection alternates between $s_1$ and $s_2$. Since we are adding a fourth string to $\alpha$ and then taking its plat closure, we may remove all beginning and trailing appearances of $\sigma_1$ or its inverse from $\alpha$. Hence we may assume that the projection of $\alpha$ is of the form $s_2(s_1 s_2)^k$ for some $k$. If $k\equiv 2$ mod 3, then it is not hard to check that the plat closure will produce a link rather than a knot. Thus we must have that $k$ is either zero or one mod 3. In either case $\alpha$ is now a Lissajous braid with frequencies $n_x=2 k+2$ and $n_y=3$. Finally, we now have that the diagram is the closure of a Lissajous arc and Theorem~\ref{closure of LA implies LP} applies.
\hfill $\square$

\section{Torus Knots}
The simplest torus knots are the $(2,2k+1)$-torus knots. Since these are 2-bridge knots, Theorem~\ref{2 bridge implies LP} implies the  following result.
\begin{corollary}
The $(2,2k+1)$-torus knot has a  Lissajous projection with $n_x=3$. 
\end{corollary}
The rest of this section will be devoted to showing that the $(3,q)$-torus knots have Lissajous projections.

We may think of the $(3,q)$-torus knot as the plat closure of the 6-string braid $\sigma_2 (\sigma_4 \sigma_3)^q \sigma_2^{-1}$. Our strategy is to focus on the five left-most strings, changing this braid into a Lissajous braid. The torus knot will then have a  projection which is the closure of a Lissajous arc and Theorem~\ref{closure of LA implies LP} will apply. 
\begin{lemma}The 6-string braids $$\sigma_2 (\sigma_4 \sigma_3)^q \sigma_2^{-1}$$ and
$$\sigma_2^{-1}\sigma_4 \sigma_1 \sigma_3^{-1}\sigma_2 \sigma_4^{-1}(\sigma_3 \sigma_2)^{q-1}\sigma_4 \sigma_3 \sigma_1 \sigma_2$$
have the same plat closure.
\label{first 6-string lemma} 
\end{lemma}
\begin{figure}
\psfrag{a}{$(\sigma_3 \sigma_2)^{q-1}$}
    \begin{center}
    \leavevmode
    \scalebox{.75}{\includegraphics{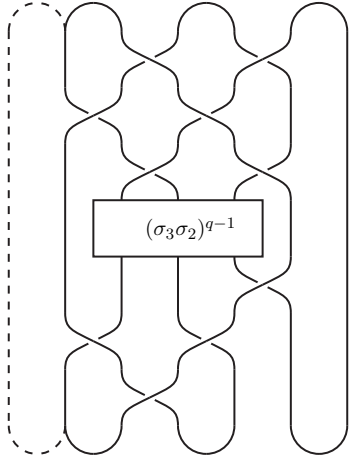}}
    \end{center}
\caption{The left-most string is picked up and moved to the right.}
\label{first 6-string lemma figure}
\end{figure}
\noindent{\bf Proof:} The proof is illustrated in Figure~\ref{first 6-string lemma figure}.
\hfill $\square$
\pagebreak
\begin{lemma} Any 5-string braid with projection 
 $$\alpha (s_3s_2)^{3n}\beta \quad \text{or} \quad \alpha (s_3s_4)^{3n}\beta$$
  can be changed to an equivalent braid with projection
$$\alpha (s_1s_3s_2s_4)^{5n}\beta.$$
\label{second 6-string lemma} 
\end{lemma}
\begin{figure}
    \begin{center}
    \leavevmode
    \scalebox{.75}{\includegraphics{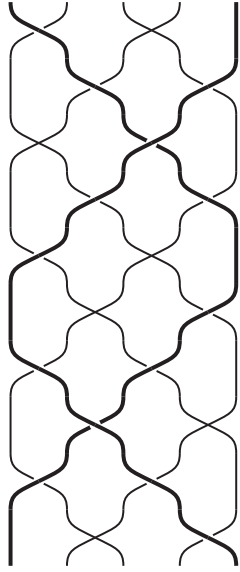}}
    \end{center}
\caption{The top-most two strings may be lifted up and pushed to the sides.}
\label{second 6-string lemma figure}
\end{figure}
\noindent{\bf Proof:} Consider a braid with projection $(s_3s_2)^3$. We may lift up the first and fifth strings and lay them down on top of the other strings as shown in Figure~\ref{second 6-string lemma figure}.

A similar proof works with $(s_3s_4)^{3n}$.
\hfill $\square$

\begin{lemma} Let $\beta$ be a 6-string braid with projection $\alpha s_3s_4s_3s_2$. Then there is a 6-string braid $\beta^\prime$ having the same plat closure as $\beta$ and having projection $\alpha s_1 s_3 s_2 s_4$.
\label{third 6-string lemma} 
\end{lemma}
\begin{figure}[h]
\psfrag{a}{$\alpha$}
    \begin{center}
    \leavevmode
    \scalebox{.75}{\includegraphics{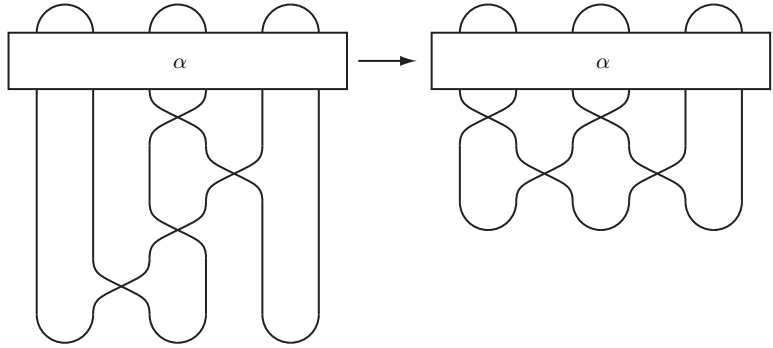}}
    \end{center}
\caption{The two lowest crossings in the left figure are moved to the left.}
\label{third 6-string lemma figure}
\end{figure}
\noindent {\bf Proof:} The proof is illustrated in Figure~\ref{third 6-string lemma figure}.

\hfill  $\square$

\begin{figure}[h]
 \psfrag{a}{$\alpha$}
    \begin{center}
    \leavevmode
    \scalebox{.75}{\includegraphics{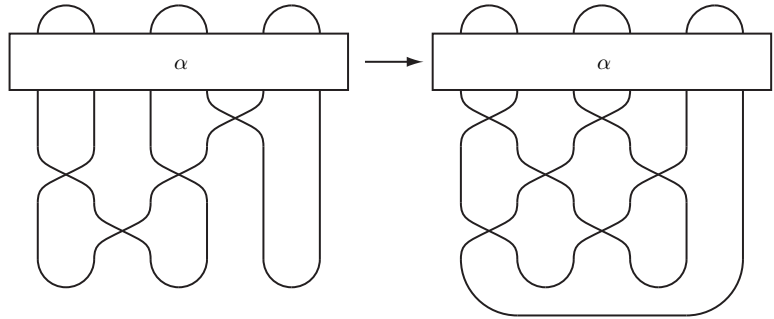}}
    \end{center}
\caption{The right-most string in the left projection is moved over (or under) the other strings as shown.}
\label{fourth 6-string lemma figure}
\end{figure}

\begin{lemma} Let $\beta$ be a 6-string braid with projection $\alpha s_4 s_3 s_1 s_2$. Then there is a 6-string braid $\beta^\prime$ with projection $\alpha s_1 s_3 s_2 s_4 s_1 s_3$  such that the plat closure of $\beta$ is equivalent to the knot obtained from $\beta^\prime$ by using ordinary plat closure at the top of $\beta^\prime$ and the modified plat closure shown in Figure~\ref{fourth 6-string lemma figure} at the bottom of $\beta^\prime$.
\label{fourth 6-string lemma} 
\end{lemma}
\noindent {\bf Proof:} Consider the right-most string in the left projection of Figure~\ref{fourth 6-string lemma figure}. Let $c$ be the crossing involving this string and having the projection $s_4$. If $c=\sigma_4$, then we may slide the right-most string over all the other strings as shown in the figure. If $c=\sigma_4^{-1}$, then we may slide the string under all the other strings.
\hfill $\square$

 We may now use Lemmas~\ref{first 6-string lemma} to \ref{fourth 6-string lemma} to prove that all $(3,q)$-torus knots have Lissajous projections. As we have already mentioned, no torus knot is Lissajous. Thus the height function corresponding to these projections cannot be given by a single cosine function.
\begin{theorem}
The $(3,q)$-torus knot has a  Lissajous projection with frequencies $5$ and $10(\frac{q-1}{3})+7$ if $q\equiv 1$ mod 3, or $5$ and $10(\frac{q-2}{3})+4$ if $q\equiv 2$ mod 3. 
\label{3 string torus knot have LP} 
\end{theorem}
\noindent{\bf Proof:} The torus knot is the plat closure of the 6-string braid $\beta=\sigma_2 (\sigma_4 \sigma_3)^q \sigma_2^{-1}$ where $q$ is relatively prime to 3.

\noindent{\bf Case I: $q=3n+1$}

By Lemma~\ref{first 6-string lemma}, we may replace $\beta$ with a braid whose projection is now
$$s_2s_4s_1s_3s_2s_4(s_3s_2)^{3n}s_4s_3s_1s_2.$$
If we now apply the first part of Lemma~\ref{second 6-string lemma} we may further alter the braid to one whose projection is 
$$s_2s_4s_1s_3s_2s_4(s_1s_3s_2s_4)^{5n}s_4s_3s_1s_2.$$
Finally, using Lemma~\ref{fourth 6-string lemma} we may change the braid so that its projection is
$$s_2s_4s_1s_3s_2s_4(s_1s_3s_2s_4)^{5n}s_1s_3s_2s_4s_1s_3=(s_2s_4s_1s_3)^{5n+3}.$$
This last change requires that we change the way the braid is closed at the bottom as described in the lemma. We have now reached a projection which is the closure of a Lissajous arc  with frequencies $5$ and $10n+7$. Hence, by Theorem~\ref{closure of LA implies LP}, the torus knot has a Lissajous projection with these same frequencies.

\noindent{\bf Case II: $q=3 n+2$}

We begin by rewriting the projection of $\beta$ as 
$$s_2s_4(s_3s_4)^{3n}s_3s_4s_3s_2.$$
Applying the second part of Lemma~\ref{second 6-string lemma} we may change the braid so that its projection is now
$$s_2s_4(s_1s_3s_2s_4)^{5n}s_3s_4s_3s_2.$$
Finally, using Lemma~\ref{third 6-string lemma}, we may change to a projection of
$$s_2s_4(s_1s_3s_2s_4)^{5n}s_1s_3s_2s_4=s_2s_4(s_1s_3s_2s_4)^{5n+1}.$$
As in Case I, we have now reached a projection which is the closure of a Lissajous arc with frequencies $5$ and $10n+4$.
\hfill $\square$

 \end{document}